\documentclass[11pt]{amsart}
\theoremstyle{plain}
\newtheorem{thm}{Theorem}[section]
\newtheorem{theorem}[thm]{Theorem}

\newtheorem{lemma}[thm]{Lemma}
\newtheorem{corollary}[thm]{Corollary}
\newtheorem{proposition}[thm]{Proposition}
\theoremstyle{definition}
\newtheorem{remark}[thm]{Remark}

\newtheorem{notation}[thm]{Notation}

\newtheorem{definition}[thm]{Definition}

\newtheorem{setup}[thm]{Set-up}
\numberwithin{equation}{section}

\newcommand{\sC}{{\mathcal C}}

\newcommand{\sF}{{\mathcal F}}

\newcommand{\sH}{{\mathcal H}}

\newcommand{\sJ}{{\mathcal J}}
\newcommand{\sK}{{\mathcal K}}

\newcommand{\sN}{{\mathcal N}}
\newcommand{\sO}{{\mathcal O}}
\newcommand{\sP}{{\mathcal P}}

\newcommand{\sR}{{\mathcal R}}

\newcommand{\sT}{{\mathcal T}}

\newcommand{\C}{{\mathbb C}}

\newcommand{\BP}{{\mathbb P}}

\newcommand{\fg}{{\mathfrak g}}
\newcommand{\fsl}{{\mathfrak s}{\mathfrak l}}

\newcommand{\fk}{{\mathfrak k}}
\newcommand{\fl}{{\mathfrak l}}
\newcommand{\fm}{{\mathfrak m}}

\newcommand{\fn}{{\mathfrak n}}
\newcommand{\fa}{{\mathfrak a}}
\newcommand{\fh}{{\mathfrak h}}

\newcommand{\fs}{{\mathfrak s}}
\newcommand{\fc}{{\mathfrak c}}

\newcommand{\ft}{{\mathfrak t}}

\def\Gr{\mathop{\rm Gr}\nolimits}

\def\Hom{\mathop{\rm Hom}\nolimits}

\title[Minimal rational curves on  symmetric varieties]{Minimal rational curves on equivariant compactifications of symmetric spaces}

\author{Jun-Muk Hwang and Qifeng Li}

\thanks{Jun-Muk Hwang was supported by the Institute for Basic Science (IBS-R032-D1). Qifeng Li was supported by the NSFC grant No. 12201348.}

\begin{document}

\begin{abstract} Let $G/H$ be a symmetric space of a complex linear algebraic group $G$ and let $X$ be a nonsingular equivariant compactification of $G/H$.  We investigate the question: when are  minimal rational curves on $X$  orbit-closures of 1-parameter subgroups of $G$? We show that this is the case if the variety of minimal rational tangents (VMRT) at a base point in $G/H \subset X$ is Gauss-nondegenerate.   Our method combines algebraic geometry of minimal rational curves with differential geometry of symmetric spaces: orbits of 1-parameter subgroups arise as holomorphic geodesics of an invariant torsion-free affine connection on $G/H$. We prove furthermore that the Gauss-nondegeneracy of VMRT holds for nonsingular equivariant compactifications of  simple algebraic groups regarded as  symmetric spaces. In this case, we also show that the  VMRT is the closure of an adjoint orbit, which generalizes a result of Brion and Fu's on wonderful compactifications  to arbitrary equivariant compactifications.
\end{abstract}

\maketitle



\section{Introduction}
Throughout, we work in the category of complex analytic varieties and holomorphic maps. Open sets refer to Euclidean topology and open sets in Zariski topology are called Zariski-open sets. A general point of an irreducible complex analytic set $X$ means a point in a  Zariski-open subset of $X$, i.e., the complement of a closed analytic subset in $X$. For a complex Lie group $G$, a $G$-variety $X$ is  a complex analytic variety $X$ equipped with a holomorphic $G$-action. We say that a subset $S \subset X$ is $G$-stable if $g\cdot S = S$ for any $g \in G$.

 \medskip
Minimal rational curves and the associated varieties of minimal rational tangents (to be abbreviated as VMRT) have been used in solving many problems on uniruled projective varieties (see  Definition \ref{d.vmrt} and the introductory surveys \cite{Hw01} and \cite{HM99}).  Here, we study minimal rational curves on nonsingular equivariant compactifications of symmetric spaces in the following sense.

\begin{definition}\label{d.symmetric} Let $G$ be a connected linear algebraic group.
\begin{itemize} \item[(i)] A {\em symmetric space}  is a triple $(G, H, \sigma)$ where $\sigma$ is  an involutive automorphism of $G$ and $H$ is an algebraic subgroup that is contained in the fixed locus $G_{\sigma}$ and contains the identity component of $G_{\sigma}$. By abuse of terminology, we often call the associated homogeneous variety $G/H$ a symmetric space. \item[(ii)] A {\em symmetric subspace} of a symmetric space $(G, H, \sigma)$ is a triple $(G', H', \sigma')$ where $G'$ is a $\sigma$-stable subgroup of $G$ with $H'= G' \cap H$ and $\sigma' = \sigma|_{G'}$. \item[(iii)]
A projective variety $X$ equipped with a $G$-action is an {\em equivariant compactification} of the symmetric space $(G, H, \sigma)$  if the $G$-orbit of a point $o \in X$ is an open subset of $X$ which is  isomorphic to $G/H$ as $G$-varieties.  \item[(iv)] A {\em base point } $o \in X$ of an equivariant compactification $X$ means   a choice of a point $o \in X$ lying on the open $G$-orbit. \end{itemize}
\end{definition}

Note that we allow the group $G$ to be non-reductive,  following the convention of \cite{KN2}. As the most interesting cases are when $G$ is reductive, many authors (for example,  \cite{BKP} and \cite{Ti}) require  $G$ to be reductive in the definition of symmetric spaces. However, even when one is mainly interested in symmetric spaces of  reductive groups, it is convenient to consider  symmetric spaces of  non-reductive groups, as they arise naturally as symmetric subspaces (or totally geodesic submanifold). See, for example,  Proposition \ref{p.vmrt} below.

When $G$ is reductive, the quotient space $G/H$ is isomorphic to the product of the quotient spaces of the  following three types of symmetric spaces (see \cite[Section 26.5]{Ti}).
\begin{itemize}
\item[({\bf T})] Torus type: $G$ is a torus and $H$ is trivial;
\item[({\bf G})] Group type: $G = S \times S$ for a simple algebraic group $S$ and $H$  is  the diagonal subgroup ${\rm diag} S \subset S \times S$;
\item[({\bf S})] Simple type: $G$ is a simple algebraic group. \end{itemize}


Equivariant compactifications of symmetric spaces of  type ({\bf T})  are toric varieties.
Minimal rational curves on nonsingular projective toric varieties are described in  \cite[Corollary 2.5]{CFH} as follows.

\begin{theorem}\label{t.CFH} Let $X$ be a nonsingular projective toric variety. Then there are finitely many nonsingular  subvarieties $X_1, \ldots, X_k$ through a base point $o \in X$ such that $X_i $ is isomorphic to $\BP^{d_i}$ for a suitable integer $d_i$ for each $1 \leq i \leq k$ and minimal rational curves on $X$ are lines on $X_i, 1 \leq i \leq k$.
\end{theorem}

A large class of symmetric spaces of type ({\bf G}) or ({\bf S}) admit a particularly nice kind of nonsingular equivariant compactifications, called  {\em wonderful compactifications}, as discovered in \cite{DP}.
In \cite[Theorem 1.1]{BF} and \cite[Table 1]{BKP}, minimal rational curves and VMRT on wonderful compactifications  have been determined. They exhibit the following common features.

\begin{theorem}\label{t.BKP}
Let $X$ be a wonderful compactification of a symmetric space $(G, H, \sigma)$ of type $({\bf G})$ or $({\bf S})$. Let $\sK$ be a family of minimal rational curves on $X$ and let $o \in X$ be a base point.  \item[(i)] Each  member of $\sK$ passing through  $o \in X$ is the closure of an orbit of a 1-parameter subgroup of $G$.
    \item[(ii)] Each irreducible component of VMRT $\sC_o \subset \BP T_o X$ of $\sK$ at $o$ is the closure of an orbit of the isotropy action of the identity component of $H$ on $\BP T_o X$.
   \end{theorem}

To be precise, the assertions in Theorem \ref{t.BKP} have not been stated explicitly in \cite{BKP}, but can be deduced from the description of minimal rational curves there.
Do the conclusions of Theorem \ref{t.BKP} hold  for an arbitrary nonsingular equivariant compactification of any symmetric space of type ({\bf G}) or ({\bf S})?
The arguments in \cite{BF} and \cite{BKP} use the algebraic geometry  of the wonderful compactification,
involving the boundary data. It does not seem to be applicable to arbitrary compactifications.
Our goal is to use tools from differential geometry of  symmetric spaces combined with the theory of VMRT-structures to study this question. Our first result is applicable to any symmetric space:

\begin{theorem}\label{t.main1}
Let $X$ be a nonsingular equivariant compactification of a symmetric space $(G, H, \sigma)$ and let $\sK$ be a family of minimal rational curves on $X$. Assume that the VMRT $\sC_o \subset \BP T_o X$ at a base point $o$ is Gauss-nondegenerate, namely, the Gauss map of the projective variety $\sC_o$ is generically finite over its image. Then
 each  member of $\sK$ passing through  $o \in X$ is  the closure of a (holomorphic) geodesic of the canonical torsion-free affine connection on the symmetric space. In particular, it is the closure of an orbit of a 1-parameter subgroup of $G$.
\end{theorem}

       The condition that $\sC_o$ is Gauss-nondegenerate is necessary. For example, in Theorem \ref{t.CFH}, if $d_i \geq 2$,  then the VMRT  $\sC_o$  is a linear subspace and Gauss-degenerate. In this case,   most lines on $\BP^{d_i}$ are not orbit-closures of 1-parameter subgroups of the torus group.

       From Theorem \ref{t.main1}, the essential problem is to check when the variety of minimal rational tangents is Gauss-nondegenerate.
       Our next result settles this for symmetric spaces of type ({\bf G}):

\begin{theorem}\label{t.main2}
Let $X$ be a nonsingular equivariant compactification of a symmetric space of type $({\bf G})$. Assume that $X$ is not biregular to projective space $\BP^n, n = \dim H$.  Then for any family of minimal rational curves $\sK$ on $X$, \begin{itemize} \item[(i)] the VMRT $\sC_o \subset \BP T_o X$ is Gauss-nondegenerate; \item[(ii)] each member of $\sK$ through  $o$ is the closure of a 1-parameter subgroup of $G$;  and \item[(iii)]  $\sC_o \subset \BP T_o X$ is the closure of an orbit of the isotropy action of $H$ on $\BP T_o X$.  \end{itemize} \end{theorem}

As a consequence of Theorem \ref{t.main2}, we can see that $\sC_o$ is an orbit closure in $\BP \fs$ whose normalization is smooth, where $\fs$ is the Lie algebra of the simple algebraic group $S$. From this perspective, it would be interesting to classify  adjoint-orbit closures in $\BP \fs$ whose normalizations are nonsingular. Such a classification is  well-known  for nilpotent orbit closures (\cite[Proposition 5.2]{Be}), but not much is known   for  orbit closures of semisimple elements.

It would be also interesting to ask which of these orbits can be realized as the VMRT of a nonsingular equivariant compactification. The result of \cite{BF} shows that this is the case for the projectivization of the minimal nilpotent orbit of a simple Lie algebra. One curious case is the exceptional group of type ${\rm G}_2$: there is a $7$-dimensional nilpotent orbit in $\BP \fg_2$ which is singular, but its normalization is nonsingular. Can it be realized as the variety of minimal rational tangents of a nonsingular equivariant compactification of an algebraic group of type ${\rm G}_2$?

We expect that an analog of Theorem \ref{t.main2} holds   for symmetric spaces of type ({\bf S}).  But  extending our arguments for Theorem \ref{t.main2} to symmetric spaces of type $({\bf S})$ would require substantial new ideas. We leave it for future investigation.

Let us give an outline of the paper.  We need some background material on differential geometry to prove Theorem \ref{t.main1}, which we review  in Section \ref{s.DG}. It consists of three parts: the theory of cone structures (\ref{ss.cone}), the theory of torsion-free affine connections (\ref{ss.affine}) and the theory of symmetric spaces (\ref{ss.symmetric}). Using these ideas, we give the proof of Theorem \ref{t.main1} in Section \ref{s.vmrt}. The proof of Theorem \ref{t.main2} is given in Section \ref{s.group}, modulo a technical lemma, Lemma \ref{l.technical}, on the geometry of the adjoint action of simple algebraic groups. The proof of this  technical lemma is rather involved and is given in Section \ref{s.adjoint}.

\section{Review of  results from differential geometry}\label{s.DG}
\subsection{Cone structures and conic connections}\label{ss.cone}

\begin{notation}\label{n.cone}
For a vector space $V$, its projectivization $\BP V$ is the set of 1-dimensional subspaces of $V$. For a projective subvariety $Z \subset \BP V$, its affine cone is denoted by $\widehat{Z} \subset V$. For a point $\alpha \in Z$, the affine tangent space of $\widehat{Z}$ at a nonzero point  $v \in \widehat{\alpha}$ is denoted by $T_{\alpha} \widehat{Z} = T_v \widehat{Z} \subset V.$ \end{notation}

\begin{definition}\label{d.cone}
A {\em cone structure} on a complex manifold $M$ is a  closed  analytic subset $\sC \subset \BP TM$ such that all irreducible components of $\sC$ have the same dimension and the natural projection $\pi: \sC \to M$ restricted to a dense Zariski-open subset $\sC^{\rm sm} \subset \sC$ is a submersion. \begin{itemize} \item[(i)] For each $x \in X$, we denote by $\sC_x \subset \BP T_x M$ the projective subvariety $\pi^{-1}(x)$ and by $\sC_x^{\rm sm}$ the intersection $\sC_x \cap \sC^{\rm sm}$.  We have natural vector subbundles $\sT \subset \sJ \subset \sP \subset T\sC^{\rm sm}$ whose fibers at $\alpha \in \sC^{\rm sm}_x$ are defined by \begin{eqnarray*} \sT_{\alpha} &=& ({\rm d}_{\alpha} \pi)^{-1}(0) \\
\sJ_{\alpha} &=& ({\rm d}_{\alpha} \pi)^{-1}(\widehat{\alpha}) \\
\sP_{\alpha} &=& ({\rm d}_{\alpha} \pi)^{-1}(T_{\alpha} \widehat{\sC}_x), \end{eqnarray*}
in terms of  the differential ${\rm d}_{\alpha} \pi: T_{\alpha} \sC^{\rm sm} \to T_x M$  of the submersion $\pi: \sC^{\rm sm} \to M$. \item[(ii)]  A line subbundle $\sF \subset \sJ$ is called a {\em conic connection} on $\sC$ if $\sJ = \sT \oplus \sF$. \item[(iii)] The {\em Cauchy characteristic}  of $\sP$ is a vector subbundle ${\rm Ch}(\sP)$  of $\sP$ on a Zariski-open subset $\sC^{\rm Cauchy}$ of $\sC^{\rm sm}$ such that a  holomorphic vector field $v$ on an open subset $O \subset \sC^{\rm Cauchy}$ is a local section of ${\rm Ch}(\sP)|_O$ if  the Lie bracket $[v, w]$ is a  section of $\sP|_O$ for any holomorphic section $w$ of $\sP|_O$.  It is easy to see that the distribution
${\rm Ch}(\sP)$ on $\sC^{\rm Cauchy}$ is integrable, defining a foliation on $\sC^{\rm Cauchy}$.
\item[(iv)] A conic connection $\sF$ is called a {\em characteristic conic connection} if $\sF \subset {\rm Ch}(\sP)$ when restricted to $\sC^{\rm Cauchy}$. \end{itemize} \end{definition}

The following is  from \cite[Proposition 1]{HM04}.

\begin{lemma}\label{l.easy}
For any conic connection $\sF$ on a cone structure $\sC \subset \BP TM$,  the sheaves of local holomorphic sections on $\sC^{\rm sm}$ satisfy
$$[\sO(\sT), \sO(\sF)] = \sO(\sP),$$  where the bracket stands for Lie brackets of local holomorphic vector fields on $\sC^{\rm sm}$.
\end{lemma}

\begin{definition}\label{d.Gauss} Let $\sC \subset \BP TM$ be a cone structure on a complex manifold $M$. Write $p := \dim \sC - \dim M $ for the dimension of the fiber $\sC_x$ at a general point $x \in X$ and let ${\rm Gr}(p+1, TM)$ be the Grassmannian bundle of $(p+1)$-dimensional subspaces in $T M$.   The {\em Gauss map} $\gamma: \sC^{\rm sm} \to {\rm Gr}(p+1, T M)$ is a morphism  sending $\alpha \in \sC^{\rm sm}_x$ to  $$\gamma(\alpha) = [ T_{\alpha} \widehat{\sC}_x] \in {\rm Gr}(p+1, T_x M).$$
\end{definition}

In the following proposition, (i) is from  \cite[Proposition 5]{HM04} and (ii) is from \cite[Proposition 4]{HM04}.

\begin{proposition}\label{p.HM}
Assume that a cone structure $\sC \subset \BP TM$ has
 a characteristic conic connection $\sF$. Then the following holds.
\begin{itemize} \item[(i)] For a germ $A$ of a leaf of the foliation ${\rm Ch}(\sP)$ on $ \sC^{\rm Cauchy}$, its projection $\pi (A)$ is a germ of a submanifold in $M$ such that the germ $A$ coincides with the germ of  $ \BP T(\pi (A)) \subset \BP T M$ as submanifolds in $\BP TM$.
\item[(ii)]  If the Gauss map $\gamma$ in Definition \ref{d.Gauss} is generically finite over its image, then $\sF = {\rm Ch}(\sP)$ on $\sC^{\rm Cauchy}$. Consequently, if the dimension of general fibers of the Gauss map $\gamma$   is zero, then  $\sF$ is the only characteristic conic connection on $\sC$. \end{itemize}
\end{proposition}

\subsection{Isotrivial cone structures and affine connections}\label{ss.affine}

In this section, we assume basic knowledge of affine connections.
The standard reference is \cite[Chapters II and III]{KN1} (see also \cite[Sections 2 and 3]{HL24}).    Although the results in \cite{KN1} are formulated for $C^{\infty}$-connections on $C^{\infty}$-manifolds over real numbers, it is straightforward to adapt any local results there to the setting of holomorphic connections on complex manifolds. By abuse of terminology,  an affine connection $\nabla$ on a complex manifold $M$ means  a linear connection arising from a principal connection on the frame bundle of $M$ and also the associated covariant derivative operator.

\begin{definition}\label{d.isotrivial}
Fix a vector space $V$ of dimension $n$ and let $M$ be a complex manifold of dimension $n$. For a projective variety $ Z \subset \BP V$, a cone structure $\sC \subset \BP TM$ is {\em $Z$-isotrivial} if the fiber $\sC_x \subset \BP T_x M$ at any $x \in M$ is isomorphic to $Z \subset \BP V$ by a linear isomorphism $T_x M \cong V$. A cone structure is {\em isotrivial} if it is $Z$-isotrivial for some $Z \subset \BP V$. \end{definition}

\begin{definition}\label{d.frame}
Let $\sC \subset \BP TM$ be a $Z$-isotrivial cone structure from Definition \ref{d.isotrivial}. Let $${\rm Aut}(\widehat{Z})  := \{ g \in {\rm GL}(V) \mid g \cdot \widehat{Z} = \widehat{Z}\} \ \subset \ {\rm GL}(V)$$ be the linear automorphism group of the affine cone $\widehat{Z} \subset V$.
\begin{itemize} \item[(i)] For each $x \in M$, let ${\bf F}^{\sC}_x$ be the set of linear isomorphisms from $V$ to $T_x M$ that sends $\widehat{Z}$ to $\widehat{\sC}_x$.
The union ${\bf F}^{\sC}:= \cup_{x \in M} {\bf F}^{\sC}_x$ with the natural projection $\varpi: {\bf F}^{\sC} \to M$ is an ${\rm Aut}(\widehat{Z})$-principal bundle on $M$, called the {\em bundle of $\sC$-frames}.
\item[(ii)] A {\em principal connection} on ${\bf F}^{\sC}$ is a vector subbundle $\sH \subset T ({\bf F}^{\sC})$ invariant under the right ${\rm Aut}(\widehat{Z})$-action on ${\bf F}^{\sC}$ such that $T({\bf F}^{\sC})= {\rm Ker} ({\rm d} \varpi) \oplus \sH$. \end{itemize} \end{definition}

\begin{remark}
Note that ${\bf F}^{\sC}$ in Definition \ref{d.frame} is called in \cite[Definition 3.2]{HL24} the {\em associated ${\rm G}$-structure} of the isotrivial cone structure, with the structure group $ {\rm Aut}(\widehat{Z})$. Since we do not need the theory of G-structures here, we decide to use the simpler name, the bundle of $\sC$-frames. \end{remark}

\begin{proposition}\label{p.connection}
In the setting of Definition \ref{d.frame}, we have the following. \begin{itemize}
\item[(i)] A principal connection $\sH$ on ${\bf F}^{\sC}$ induces an affine connection $\nabla^{\sH}$ on $M$ and  a conic connection $\sF^{\sH}$ on $\sC \subset \BP TM$ such that leaves of $\sF^{\sH}$ descend to (holomorphic) $\nabla^{\sH}$-geodesics on $M$.
\item[(ii)] An affine connection $\nabla$ on $M$ is induced by a principal connection $\sH$ on ${\bf F}^{\sC}$, namely, it is of the form $\nabla^{\sH}$ in (i), if and only if $\sC \subset \BP TM$ is preserved under the $\nabla$-parallel displacements.
    \item[(iii)] If the affine connection $\nabla^{\sH}$ in (i) is torsion-free, namely,
     $$\nabla^{\sH}_{\vec{u}} \vec{v} - \nabla^{\sH}_{\vec{v}} \vec{u} \ = \ [ \vec{u}, \vec{v}]$$
     for all local holomorphic vector fields $\vec{u}, \vec{v}$ on $M$,
     then $\sF^{\sH}$ is a characteristic conic connection.
        \item[(iv)] Assume that $\sF^{\sH}$ in (i) is a characteristic conic connection. Then for a germ $A$ of a leaf of the foliation ${\rm Ch}(\sP),$ the submanifold germ $\pi(A) \subset M$ is totally geodesic with respect to $\nabla^{\sH}$, namely, for any tangent vector $v$ to $\pi(A)$, the $\nabla^{\sH}$-geodesic in the direction of $v$ lies on the submanifold $\pi(A)$. \end{itemize}
        \end{proposition}

        \begin{proof} In (i), that a principal connection induces an affine connection is standard (for example, \cite[Chapter III]{KN1}) and that the geodesics of $\nabla^{\sH}$ determines a conic connection is from \cite[Proposition 3.4]{HL24}. (ii) is a consequence of \cite[Theorem 7.1 in Chapter II]{KN1}. (iii) is a special case of \cite[Theorem 3.8]{HL24}, when the torsion $\tau^{\sH}$ is identically zero. In (iv), the assumption says  $\sF^{\sH} \subset {\rm Ch}(\sP).$ By (i), leaves of $\sF^{\sH}$ in $A$ descend to $\nabla^{\sH}$-geodesics to cover $\pi(A)$. Moreover, such $\nabla^{\sH}$-geodesics span all tangent directions of $\pi(A)$ by Proposition \ref{p.HM} (i). This proves (iv). \end{proof}

\begin{remark} In  \cite[Definition 3.1]{HL24}, cone structures $\sC \subset \BP TM$ are assumed be smooth, while our cone structure can have singularity. However, the results from \cite{HL24} used in the proof of Proposition \ref{p.connection} are concerned only with what is happening on  a Zariski-open subset  of $\sC$, so we may apply it to $\sC^{\rm sm} \subset \sC$. \end{remark}

\subsection{Torsion-free affine connections on symmetric spaces}\label{ss.symmetric}
We need to recall some properties of the canonical torsion-free affine connections on symmetric spaces. In the next proposition,  (i) and (ii) are from  \cite[Theorem 3.2 in Chapter XI]{KN2}  and (iii) is from \cite[Theorem 4.3 in Chapter XI]{KN2}.

\begin{proposition}\label{p.nabla}
Let $(G, H, \sigma)$ be a symmetric space and  $M= G/H$ be the associated homogeneous complex manifold with the base point $o=[H] \in G/H$. Let $\fg = \fh + \fm$ be the $\sigma$-eigenspace decomposition of the Lie algebra with a natural identification $\fm = T_o M$. Then there exists a unique $G$-invariant torsion-free affine connection $\nabla$ on $M$ with the following properties. \begin{itemize} \item[(i)]
For $v \in \fm = T_o M$,  the orbit of the one-parameter subgroup $t \in \C \mapsto \exp(t v) \in G$ through $o \in M$ is the $\nabla$-geodesic through $o$ in the direction of $v$.
\item[(ii)] In (i), let  $\mathbf{e}^{t v}: M \to M$ be the automorphism of $M$ given by the action of $\exp(tv)   \in G$ inducing the differential $${\rm d}_o \mathbf{e}^{tv}: T_o M \to T_{e^{tv}(o)}M.$$ Then the parallel translate of a vector $w \in T_o M$ with respect to $\nabla$ along the geodesic in (i) to the point $\mathbf{e}^{tv}( o)$ coincides with ${\rm d}_o
\mathbf{e}^{tv} (w)$.
\item[(iii)] Let $o \in M' \subset M$ be a closed totally geodesic submanifold. Then there exists a vector subspace $\fm' \subset \fm$ satisfying $[[\fm', \fm'], \fm'] \subset \fm'$ such that $M'$ is the symmetric subspace of $M$ determined by the Lie subalgebra $\fg' \subset \fg$ with the involution $\sigma'$ defined by  $$\fg':= \fh' + \fm', \ \fh': = [ \fm', \fm'] \mbox{ and } \sigma' := \sigma|_{\fg'}.$$ \end{itemize} \end{proposition}

It has the following consequence.

\begin{proposition}\label{p.cone}
For a symmetric space $(G, H, \sigma)$, the tangent bundle $TM$ of $M = G/H$ can be identified with the vector bundle $G \times_H \fm$ associated with the $H$-principal bundle $G \to G/H$ via the isotropy action of $H$ on the $(-1)$-eigenspace  $\fm \subset \fg$ of $\sigma.$ Let $\nabla$ be the canonical torsion-free affine connection on $M = G/H$ from  Proposition \ref{p.nabla}. \begin{itemize}
\item[(i)] For any $H$-stable projective subvariety $Z \subset \BP \fm$, the fiber subbundle $$G \times_H  Z \subset G \times_H \BP \fm = \BP T M,$$ is a $Z$-isotrivial cone structure on $M$, which is $G$-stable under the natural $G$-action on $\BP TM$.
    \item[(ii)] A $G$-stable cone structure $\sC \subset \BP TM$ on $M$ is $Z$-isotrivial for $Z := \sC_o \subset \BP \fm$, the fiber at a base point $o \in M,$ and the fiber subbundle $$G \times_H \sC_o \subset G \times_H \BP \fm = \BP T M$$ coincides with $\sC \subset \BP TM.$
    \item[(iii)] For a $G$-stable cone structure $\sC \subset \BP TM,$ the canonical connection $\nabla$ agrees with the affine connection $\nabla^{\sH}$ arising from a principal connection $\sH$ on the bundle ${\bf F}^{\sC}$ of $\sC$-frames.   \item[(iv)] The $\nabla$-geodesics induce a characteristic conic connection on any $G$-stable cone structure $\sC \subset \BP TM$. \end{itemize} \end{proposition}

        \begin{proof}
     (i) is immediate.
        As $G$ acts transitively on $M$, two $G$-stable subsets in $\BP TM$ coincide if and only if their fibers coincide  at a base point $o \in M,$ proving (ii).   Proposition \ref{p.nabla} (ii) implies that $\nabla$-parallel displacements preserve any $G$-stable cone structure $\sC \subset \BP TM$. Thus   $\nabla$ comes from a principal connection on ${\bf F}^{\sC}$ by Proposition \ref{p.connection} (ii), proving (iii).   (iv) follows from (iii) and Proposition \ref{p.connection} (iii), because $\nabla$ is torsion-free.  \end{proof}

\section{VMRT-structures of equivariant compactifications of symmetric spaces}\label{s.vmrt}
Many interesting examples of cone structures with characteristic conic connections arise from minimal rational curves as follows.

\begin{definition}\label{d.vmrt}
Let $X$ be a nonsingular projective variety.
\begin{itemize} \item[(i)] An irreducible component $\sK$ of the space ${\rm RatCurves}(X)$ of rational curves on $X$ is called a {\em family of minimal rational curves} on $X$ if for the associated universal family morphisms $$\sK \ \stackrel{\rho}{\leftarrow} \ {\rm Univ}_{\sK} \ \stackrel{\mu}{\rightarrow} \ X,$$ the fiber $\mu^{-1}(x)$ at a general point $x \in X$  is nonempty and projective. \item[(ii)] Fix a family $\sK$ of minimal rational curves on $X$. Then there exists a nonempty Zariski-open subset $M \subset X$ with a cone structure
 $\sC \subset \BP TX$, called the {\em VMRT-structure} associated to the family $\sK$, such that for a general point $x \in M$, a general point of $\sC_x$ is  the tangent direction of a member of $\sK$ through $x$ and any member of $\sK$ through $x$ which is smooth at $x$ has its tangent direction belonging to $\sC_x$.
We call the projective variety (with finitely many irreducible components) $\sC_x \subset \BP T_x X$ the {\em VMRT} of $\sK$ at $x \in X$.
\end{itemize} \end{definition}

The following is from \cite[Corollary 1]{HM04}.

\begin{proposition}\label{p.birational}
In Definition \ref{d.vmrt}, the normalization of the VMRT $\sC_x $ at a general point $x \in X$ is nonsingular. \end{proposition}

The following is from  \cite[Proposition 8]{HM04}.

\begin{proposition}\label{p.vmrtconnection}
Let $\sC \subset \BP TM$ be as in  Definition \ref{d.vmrt}, the VMRT structure associated to a family of minimal rational curves on a nonsingular projective variety $X.$    Then the tangent vectors of general members of $\sK$ determine
 a characteristic conic connection $\sF^{\sK} \subset T \sC^{\rm sm}$ of the cone structure $\sC \subset \BP TM.$  \end{proposition}

The following is from  \cite[Propositions 9 and 12]{HM04}

\begin{proposition}\label{p.Cauchy}
Let $\sC \subset \BP TM$ be the VMRT structure in Definition \ref{d.vmrt} and let $p:= \dim \sC - \dim M$ be the fiber dimension of $\sC \to M$.  Let $A$ be a general leaf of ${\rm Ch}(\sP)$ on $\sC^{\rm Cauchy}$.  Then the closure $\bar{S} \subset X$ of the image $S:= \pi(A) \subset X$ is a projective subvariety of $X$ such that
\begin{itemize}
\item[(i)]  at a general point $x \in \bar{S}$, the linear subvariety $\BP T_x \bar{S} \subset \BP T_x X$ is the closure of a fiber of the Gauss map $\gamma_x: \sC^{\rm sm}_x \to {\rm Gr}(p+1, T_x X)$; and
\item[(ii)]
there exists a finite morphism $\BP^d \to \bar{S}$ where $d = \dim \pi(A),$ which sends lines on $\BP^d$ to members of $\sK$ lying on $\bar{S}$. \end{itemize}
\end{proposition}

When a connected group $G$ acts on a nonsingular projective variety  $X,$ it preserves each  irreducible component of ${\rm RatCurves}(X)$, inducing a natural $G$-action on a family $\sK$ of minimal rational curves and the associated VMRT-structure $\sC \subset \BP TM.$ Hence the following is a direct corollary of  Proposition \ref{p.cone}.

\begin{proposition}\label{p.G-stable}
Let $X$ be a nonsingular equivariant compactification of a symmetric space $M = G/H.$ We choose a base point $o \in X$ identifying  $M$ with the $G$-orbit of $o \in X$. Let $\sC \subset \BP TM$ be the VMRT-structure given by a family of minimal rational curves on $X$. Then the fiber $\sC_o \subset \BP T_o M$ is $H$-stable and
the restriction $\sC|_M$ is a $G$-stable cone structure. \end{proposition}

We have the following results on minimal rational curves of symmetric varieties, which proves Theorem \ref{t.main1}.

\begin{proposition}\label{p.vmrt}
Let $X$ be a  nonsingular equivariant compactification of a symmetric space $(G, H, \sigma)$. Let $o \in X$ be a  point whose $G$-orbit is an open subset $M \subset X$.   Let $\sC \subset \BP TX$ be the VMRT-structure given by a family $\sK$ of minimal rational curves on $X$.
\begin{itemize} \item[(i)] For the subvariety $\bar{S} \subset X$ through the base point $o \in M \subset X$ in Proposition \ref{p.Cauchy}, the open subset $\bar{S} \cap M$ is a closed totally geodesic submanifold of the symmetric space $M$. In particular, in terms of the Lie algebra $\fg = \fh + \fm$ of $(G, H, \sigma)$, the tangent space $T_o S \subset \fm$ satisfies $[[T_o S, T_o S], T_o S] \subset T_o S$ in $\fg$.
\item[(ii)] If the Gauss map $\gamma$ of $\sC_o$ is generically finite over its image, then each member of $\sK$ through $o \in M$ is the closure of a $\nabla$-geodesic for the canonical torsion-free connection $\nabla$ in Proposition \ref{p.nabla}, hence the closure of an  orbit of a 1-parameter subgroup $\exp (\C v)$ for some $ v \in \fm$.
  \end{itemize}
\end{proposition}

\begin{proof}
By Proposition \ref{p.G-stable}, the VMRT-structure $\sC \subset \BP TM$  is a $G$-stable cone structure. Thus, by Proposition \ref{p.cone} (iii), the canonical torsion-free connection $\nabla$ agrees with $\nabla^{\sH}$ for a principal connection $\sH$ on the bundle ${\bf F}^{\sC}$ of $\sC$-frames. By Proposition \ref{p.connection} (iii), the conic connection $\sF^{\sH}$ is a characteristic conic connection.  Then
 Proposition \ref{p.connection} (iv) implies that the submanifold $\bar{S} \cap M$ in Proposition \ref{p.Cauchy} is totally geodesic with respect to $\nabla$, proving
 the first statement in (i).  The second statement in (i) follows Proposition \ref{p.nabla} (iii).

 If the Gauss map of $\sC_o$ is generically finite over its image, the uniqueness of the characteristic conic connection in Proposition \ref{p.HM} implies that the characteristic conic connection $\sF^{\sK}$ in Proposition \ref{p.vmrtconnection} coincides with the conic connection $\sF^{\sH}$ in Proposition \ref{p.connection} (i), which is a characteristic conic connection by Proposition \ref{p.connection} (iii). It follows that the germs of members of $\sK$ through $o$ are $\nabla$-geodesics, which are closures of orbits of 1-parameter subgroups by Proposition \ref{p.nabla} (i).
 \end{proof}

\section{Minimal rational curves on a nonsingular equivariant compactification of a simple algebraic group}\label{s.group}

In the rest of the paper, we change  our notational convention as follows. We fix  a simple algebraic group $G$. Let  $\fg$ be the Lie algebra of $G$ and regard $\BP \fg$ as a $G$-variety by the adjoint action. Let $\sN \subset \fg$ be the cone of nilpotent elements of $\fg$. It is well-known that there are only finitely many $G$-orbits in $\sN$.

We view $G$ as a symmetric space of the group $G \times G$ by the left and right translations with the involution $\sigma$ on $G \times G$ given by $\sigma (g_1, g_2) = (g_2, g_1).$ We consider an equivariant compactification $X$ of $G$ as a symmetric space. In this case, there exists a unique family $\sK$ of minimal rational curves on $X$ by  \cite[Proposition 2.6]{BF} and consequently, a uniquely determined VMRT-structure $\sC \subset \BP TM, M = (G \times G)/{\rm diag} G \cong G$.
Our main result is the following.

 \begin{theorem}\label{t.group}
 Let $X$ be an equivariant compactification of the symmetric space $(G \times G, {\rm diag} G, \sigma),$ not biregular to $\BP^{\dim G}$.  Then  the VMRT $\sC_o \subset \BP T_o X$ at a base point $o \in X$ is Gauss-nondegenerate.
  \end{theorem}

This is exactly Theorem \ref{t.main2} (i). We have   the following  corollary of Theorem \ref{t.group}, which is exactly Theorem \ref{t.main2} (ii) and (iii).

 \begin{corollary}\label{c.group}
In Theorem \ref{t.group}, let $o \in X$ be a base point.
\begin{itemize} \item[(i)] All minimal rational curves through  $o \in X$ is the closure of the orbit of a 1-parameter subgroup of $G$.
 \item[(ii)] Identify $T_o X$ with the Lie algebra $\fg$ of $G$. Then the VMRT $\sC_o \subset \BP \fg$ is the closure of a $G$-orbit in $\BP \fg$. \end{itemize}
     \end{corollary}

 \begin{proof}[Proof of Corollary \ref{c.group}]
(i) is a direct consequence of Theorem \ref{t.main1} and Theorem \ref{t.group}.

For (ii), let  $\widehat{\sC}_o \subset \fg$ be  the affine cone of the VMRT at $o$ and let $Z$ be an irreducible component of $\widehat{\sC}_o$, which should be stable under the adjoint action of $G$.
 If  $Z \subset \sN$, then  $Z$ must be the closure of one of the finitely many $G$-orbits  in $\sN$.

If $Z \not\subset \sN$, then for a general element $ v \in Z$,  the 1-parameter group $\exp(\C v) \subset G$ is  a 1-dimensional algebraic torus, because its closure in $X$ is an algebraic curve.  Since there are only countably many 1-dimensional tori in $G$ up to conjugation, we see that $Z$ is the closure of the adjoint orbit of $\C v$.

Since $\widehat{\sC}$ is irreducible and $\widehat{\sC}_o$ is of pure dimension, $(G\times G)\cdot Z$ is the restriction of $\widehat{\sC}$ to the open orbit in $X$. But the fiber of $(G\times G)\cdot Z$ is exactly $Z$, so we have $\widehat{\sC}_o=Z$.
\end{proof}

To prove Theorem \ref{t.group},  we need the following three lemmata.

\begin{lemma}\label{l.nilpotent}
Let $Z \subset \BP \fg$ be a $G$-stable projective subvariety satisfying $\widehat{Z} \subset \sN$.    \begin{itemize}
      \item[(i)] If the normalization of $Z$ is nonsingular, then
          \begin{itemize} \item[(a)] either  $Z$ is the unique closed $G$-orbit in $\BP \fg$, the orbit of the root space of a long simple root; or  \item[(b)]  $\fg$ is of type ${\rm G}_2$ and  $Z$ is the closure of the 7-dimensional orbit of the root space of  a short simple root. \end{itemize}
            \item[(ii)] The projective subvariety $Z \subset \BP \fg$ in both (a) and (b) of  (i) is Gauss-nondegenerate.  \end{itemize} \end{lemma}

                \begin{proof} (i) is \cite[Proposition 5.2]{Be}.
                Since $Z$ in (a) is nonsingular and nonlinear, it is Gauss-nondegenerate (for example, by \cite[Corollary 3.4.18]{La}).

                It remains to check the Gauss-nondegeneracy of $Z$ in  (b).
                Fix a Cartan subalgebra and a system of simple roots of $\fg$ of type $G_2$ such that $\alpha_1$ is a short simple root and $\alpha_2$ is a long simple root. The positive roots of $\fg$ are $$\alpha_1, \alpha_2, \alpha_1+\alpha_2, 2 \alpha_1+\alpha_2, 3 \alpha_1+\alpha_2, 3\alpha_1+2\alpha_2.$$ Given a root $\beta$ of $\fg$, we denote by $\fg_\beta$ the corresponding root space and by $\ft_\beta:=[\fg_\beta, \fg_{-\beta}]$ the corresponding 1-dimensional subspace of the  Cartan subalgebra. Let $x \in Z$ be the point $\BP \fg_{\alpha_1}$ such that $Z$ is the closure of the $G$-orbit of $x$.
                            A direct calculation shows that the Lie algebra of the isotropic subgroup $G_x$ of the $G$-action at $x\in Z$ is
\begin{eqnarray*}
\lefteqn{ \{v\in\fg\mid [v, \fg_{\alpha_1}]\subset\fg_{\alpha_1}\}  } \\ & = &\ft_{\alpha_1}+\ft_{\alpha_2}+\fg_{\alpha_1}+\fg_{3\alpha_1+\alpha_2}+\fg_{3\alpha_1+2\alpha_2}+\fg_{-\alpha_2}+\fg_{-3\alpha_1-2\alpha_2}.
\end{eqnarray*}
On the other hand, the tangent space $ T_x\widehat{Z}$ of the affine cone $\widehat{Z}$ at $x$ is
\begin{eqnarray*}   [\fg, \fg_{\alpha_1}]  & = & \ft_{\alpha_1}+\fg_{\alpha_1}+\fg_{\alpha_1+\alpha_2}+\fg_{2\alpha_1+\alpha_2} \\ & &
+\fg_{3\alpha_1+\alpha_2}+\fg_{-\alpha_2}+\fg_{-\alpha_1-\alpha_2}
+\fg_{-2\alpha_1-\alpha_2}.
\end{eqnarray*}
From the above two equations, it is easy to check
$$
\{v\in\fg\mid [v, \fg_{\alpha_1}]\subset\fg_{\alpha_1}\}=\{v\in\fg\mid [v, T_x\widehat{Z}]\subset T_x\widehat{Z}\},$$
which implies that the Gauss map  of $Z$ has 0-dimensional fiber through the point $x \in Z$.
\end{proof}

\begin{lemma}\label{l.PSL2}
Let $H$ be a three-dimensional simple algebraic group. Let $Y$ be a normal equivariant compactification of the symmetric space $(H \times H, {\rm diag} H, \sigma)$. Assume that there exists a finite morphism $f: \BP^3 \to Y$ such that the family of curves on $Y$ given by the images of lines of $\BP^3$ under $f$ is $(H \times H)$-stable. Then $H$ is isomorphic to ${\rm PSL}_2(\C).$ \end{lemma}

\begin{proof}
It suffices to assume that $H$ is isomorphic to ${\rm SL}_2(\C)$ and derive a contradiction.
Recall that a biholomorphic map $\varphi: U \to \widetilde{U}$ between two nonempty connected open subsets $U, \widetilde{U} \subset \BP^3$ that sends germs of lines on $U$ to germs of lines on $\widetilde{U}$ can be extended to   a  biregular automorphism of $\BP^3$ (this is a version of the fundamental theorem of projective geometry and can be also seen as a special case of Cartan-Fubini type extension theorem in \cite{HM01}).
 Thus the action of the simply connected group $H \times H$ on $Y$ can be lifted to an action of $H \times H$ on $\BP^3$ such that $f$ is $H\times H$-equivariant. The  left-translation by $H$ acts simply transitively on a Zariski-open subset $Y^o \subset Y$. Thus it acts also simply transitively on $f^{-1}(Y^o) \subset \BP^3$, which implies that $f$ must be birational. As $Y$ is normal, we conclude that $f$ is biregular and $Y$ is isomorphic to $\BP^3$. But a nonsingular equivariant compactification of ${\rm SL}_2(\C)$ as a symmetric space
cannot be isomorphic to $\BP^3$ (for example, by \cite[Theorem 4 (1)]{Ru10} and \cite[Theorem 4.1]{Ru11}), a contradiction. \end{proof}

For the next lemma, recall that we call a Lie subalgebra $\fh \subset \fg$  an {\em algebraic} Lie subalgebra if it is the Lie algebra of an algebraic subgroup of $G$.

\begin{lemma}\label{l.technical}
Let $Z\subset\BP\fg$ be a  projective $G$-stable subvariety of dimension $p < \dim \BP \fg $ and let $Z^{\rm sm} \subset Z$ be its smooth locus.  Let  $\sR\subset\Gr(p+1, \fg)$ the  closure of the image of the Gauss map $\gamma: Z^{\rm sm} \to \Gr(p+1, \fg)$. Suppose that
\begin{itemize}
\item[(a)] $Z$ is Gauss-degenerate, namely, the fibers of $\gamma$ have positive dimensions;
\item[(b)] for a general $r\in\sR$, there exists an algebraic Lie subalgebra $\fh_r \not\subset \sN$ of $\fg$ such that $\BP\fh_r$ is the closure of the Gauss fiber $\gamma^{-1}(r)$; and
    \item[(c)] the normalization of $Z$ is nonsingular.
\end{itemize}
Then  the Lie algebra $\fh_r$ in (b) is isomorphic to $\fsl_2(\C)$ and the corresponding subgroup $\exp(\fh_r) \subset G$ is isomorphic to ${\rm SL}_2(\C).$
\end{lemma}

Lemma \ref{l.technical} is of technical nature, in the sense that we do not know whether a projective subvariety $Z \subset \BP \fg$ with the described properties  actually exists. Its proof is quite involved and  is given in the next section.

\begin{proof}[Proof of Theorem \ref{t.group}]
Assuming that an irreducible component $Z \subset \BP \fg$ of $\sC_o \subset \BP T_o X$ is Gauss-degenerate, we derive a contradiction. Since $Z$ is $G$-stable and the normalization of $Z$ is nonsingular by Proposition \ref{p.birational},  Lemma \ref{l.nilpotent} (ii)  implies that $Z \not\subset \sN$.

By Proposition  \ref{p.Cauchy} (i),  the closure of a general fiber of the Gauss map of $Z$ is of the form $\BP W$ for a linear subspace $W \subset \fg$. For a general element $w \in W,$ we see from $[w, w]=0$ that the action of the subgroup $\exp(\C w) \subset G$ fixes the point $[w] \in Z$. Thus the fiber $\BP W \subset Z$ is preserved by the action of $\exp(\C w),$ because the Gauss map $\gamma$ is equivariant under the action of $\exp(\C w) \subset G.$ This implies that $[w, W] \subset W$ for a general $w \in W$ and $W \subset \fg$ is a Lie subalgebra of $\fg.$ Consequently,  the symmetric subspace  of  the symmetric space $(G \times G, {\rm diag} G, \sigma)$  corresponding to a general fiber of $\gamma$ from Proposition \ref{p.vmrt} (i)   is of the form $(H \times H, {\rm diag} H, \sigma|_H)$ for a Lie subgroup $H \subset G$. We conclude that $W$ is an algebraic Lie subalgebra of $\fg.$

Thus the conditions (a), (b) and (c) of Lemma \ref{l.technical} are satisfied. It follows that $H \subset G$ is a closed subgroup isomorphic to ${\rm SL}_2(\C)$. The closure $\bar{H}$ in $X$ of the subgroup $H \subset G \cong M \subset X$  is an equivariant compactification of $(H \times H, {\rm diag} H, \sigma|_H)$ and the family of minimal rational curves of $X$ lying on $\bar{H}$ is preserved under the action of $H\times H$.  The normalization $Y$ of  $\bar{H}$ is a normal equivariant compactification of $H$ as a symmetric space. By Proposition \ref{p.Cauchy} (ii), there is a finite morphism $\BP^3 \to Y$ such that the family of the images of lines of $\BP^3$ is $(H \times H)$-stable. Then Lemma \ref{l.PSL2} implies that $H$ is isomorphic to ${\rm PSL}_2(\C)$, a contradiction.
\end{proof}

\section{Adjoint-invariant subvarieties  fibered by Lie subalgebras}\label{s.adjoint}

 This section is devoted to the proof of Lemma \ref{l.technical}.  In particular, we study the geometry of a  $G$-stable subvariety of  $\BP \fg$ in the following setting.

\begin{setup}\label{setup} Let $Z \subset \BP \fg$ be a $G$-stable projective subvariety. Assume that there exist \begin{itemize} \item[(1)] a $G$-stable nonempty Zariski-open subset $Z^o \subset Z$; \item[(2)]   a $G$-equivariant smooth morphism $\gamma: Z^o \to R^o$ onto a positive-dimensional nonsingular algebraic variety $R^o$ equipped with a $G$-action; and \item[(3)] an algebraic Lie subalgebra $\fh_r \subset \fg$ of dimension at least 2 for each $r \in R^o,$ \end{itemize} such that \begin{itemize} \item[(i)] $\gamma^{-1}(r)$ is a Zariski-open subset of the projective linear subspace $\BP \fh_r \subset \BP \fg$ ; and
\item[(ii)] $\fh_r \cap \widehat{Z}^o$ consists of semisimple elements of $\fg$ \end{itemize} for each $r \in R^o$. \end{setup}

We recall the following elementary fact.

\begin{lemma}\label{l.countable}
Let $T \subset G$ be a maximal torus of a simple algebraic group $G$. Then there are only countably many $T$-stable algebraic Lie subalgebras in $\fg$. \end{lemma}

\begin{proof}
Let $\Phi$ be the root system of $\fg$ with respect to the Cartan subalgebra $\ft \subset \fg$ corresponding to $T \subset G$.
For a $T$-stable algebraic Lie subalgebra $\fl \subset \fg$, we have the decomposition (see \cite[Chapter 6, Section 1]{OV})  \begin{equation}\label{e.decompo}  \fl = (\fl \cap \ft) \oplus (\oplus_{\alpha \in \Phi_{\fl}} \fg_{\alpha}),\end{equation} where  $\Phi_{\fl} := \{ \alpha \in \Phi \mid \fg_{\alpha} \subset \fl\}$.
It is well-known that there are only countably many algebraic subalgebras in the Lie algebra $\ft$. So in (\ref{e.decompo}), we have only countably many possibilities of $\fl \cap \ft$ and only finitely many possibilities of $\Phi_{\fl} \subset \Phi$. It follows that there are only countably many possibilities of $\fl$. \end{proof}

\begin{proposition}\label{p.transitive}
In Set-up \ref{setup},  the $G$-action on $R^o$ has only finitely many orbits. In particular, there is a nonempty open $G$-orbit in $R^o$. \end{proposition}

\begin{proof}
Fix a  Cartan subalgebra $\ft \subset \fg$ with the corresponding maximal torus $T \subset G$. Let $R^T \subset R^o$ be the closed algebraic subset consisting of $T$-fixed points.
For each $r \in R^T$, the fiber $\gamma^{-1}(r)$ is $T$-stable, hence the Lie subalgebra $\fh_r$ is a $T$-stable algebraic subalgebra of $\fg$. By Lemma \ref{l.countable}, we see that $R^T$ is a finite set.

Thus to prove the proposition, it suffices to check that for each $r \in R^o$, there exists $a \in G$ such that $a \cdot r \in R^T$. By Set-up \ref{setup} (ii), we have a nonzero semi-simple element $v \in \fh_r$.  Let $\ft' \subset \fg$ be a Cartan subalgebra containing $v$ and choose $a \in G$ such that ${\rm Ad}_a(\ft') = \ft$. The maximal torus $T'$ corresponding to  $\ft'$ fixes $v$, hence fixes $r \in R^o$. It follows that $T$ fixes $a \cdot r,$ proving the proposition. \end{proof}

\begin{corollary}\label{c.TZ}
In Proposition \ref{p.transitive},  for a general point $r\in R$ and a general element $v \in \fh_r$, we have  $T_v \widehat{Z} = \fh_r + [\fg, v].$ \end{corollary}

\begin{proof}
The tangent map $T_v\widehat{Z}\rightarrow T_r R$ has kernel $\fh_r$ and by Proposition \ref{p.transitive} this map sends $[\fg, v]$ onto $T_rR$. It follows that $T_v \widehat{Z} = \fh_r + [\fg, v].$
\end{proof}

Given a vector subspace $W \subset \fg$,  we  write \begin{eqnarray*} \fc_{\fg}(W) & := & \{ u \in \fg \mid [u, W]=0\} \\
\fn_{\fg}(W) &:= & \{ u \in \fg \mid [u, W] \subset W\} \end{eqnarray*} for the centralizer and the normalizer of $W$, respectively. Also, we write $\fc_{\fg}(s) = \fc_{\fg}(\C s)$ and $\fn_{\fg}(s) = \fn_{\fg}(\C s)$ for an element $s \in \fg$.

\begin{proposition}\label{p.TZ} In Set-up \ref{setup},  fix a general point $r \in R^o$,  a general point $s \in \fh_r$ and a Cartan subalgebra $\ft \subset \fg$ containing $s$.
 Write $$ \fh :=  \fh_r,  \ \fh_0 := \fh \cap \ft, \mbox{ and }
 \fn  :=  \fn_{\fg}(\fh),$$ to simplify the notation. By the $G$-equivariance of $\gamma$, we have  $\ft \subset \fn$ and the subalgebra $\fn$ is $\ft$-stable. Thus we can choose a $\ft$-stable subspace $\fm \subset \fg$ such that $\fg = \fn \oplus \fm$. Then $T_s \widehat{Z} = \fh \oplus \fm$ and it is a $\ft$-stable subspace of $\fg$. \end{proposition}

 \begin{proof} From Corollary \ref{c.TZ}, $$ T_s \widehat{Z} = \fh + [\fg, s] = \fh + [\fn, s] + [\fm, s] \ \subset \fh + \fm.$$ The Lie algebra of the isotropy subgroup $N \subset G$ of the point $r \in R^o$ must be the normalizer of $\fh \subset \fg$. Thus $\dim R^o = \dim \fg - \dim \fn$ and   $\dim T_s \widehat{Z} = \dim \fh + \dim R^o = \dim \fh + \dim \fm$. It follows that $T_s \widehat{Z} = \fh \oplus \fm.$ \end{proof}

 Recall the following elementary fact.

 \begin{lemma}\label{l.toral}
 For a subspace $W \subset \ft$ of the Cartan subalgebra, we have $\fn_{\fg}(W) = \fc_{\fg}(W).$ \end{lemma}

 \begin{proof}
 Choose a decomposition $\fg = W \oplus V$ into $\ft$-stable subspaces. For any $a \in \fn_{\fg}(W),$ write $a = w + v$ with $ w \in W, v \in V$. Since $[w, W]=0$ and $[a, W] \subset W$, we have $[v, W] = [a-w, W] \subset W$. But $[v, W] \subset V$ because $V$ is $\ft$-stable. Thus $[v, W]=0$ and $a \in \fc_{\fg}(W).$ \end{proof}

 \begin{proposition}\label{p.tstable}
 In Proposition \ref{p.TZ}, \begin{itemize}
 \item[(i)] $\fn_{\fg}(s) = \fc_{\fg}(s) \subset \fn.$
 \item[(ii)] The subalgebras $\fh, \fc_{\fg}(s), \fc_{\fg}(\fh_0), \fc_{\fg}(\fh)$ and $\fn$ of $\fg$ are $\ft$-stable.
     \item[(iii)] $\fc_{\fg}(\fh) + \ft \subset \fc_{\fg}(\fh_0) \subset \fc_{\fg}(s).$
     \end{itemize} \end{proposition}

     \begin{proof} The equality in (i) is from Lemma \ref{l.toral} and the inclusion in (i) is from the $G$-equivariance of $\gamma.$

     The inclusions of Lie algebras $\ft \subset \fc_{\fg}(\fh_0) \subset \fc_{\fg}(s) \subset \fn$ imply that $\fc_{\fg}(\fh_0), \fc_{\fg}(s)$ and $\fn$ are $\ft$-stable. From $[\ft, \fh] \subset [\fn, \fh] \subset \fh$, we see that $\fh$ is $\ft$-stable. Then
     \begin{eqnarray*} [[\ft, \fc_{\fg}(\fh)], \fh] & \subset & [[\ft, \fh], \fc_{\fg}(\fh)] + [\ft, [\fc_{\fg}(\fh), \fh]] \\ & \subset & [\fh, \fc_{\fg}(\fh)] = 0 \end{eqnarray*}
     shows that $[\ft, \fc_{\fg}(\fh)] \subset \fc_{\fg}(\fh).$ This proves (ii).

     (iii) is immediate from $s \in \fh_0 = \fh \cap \ft$ and $[\ft, \ft]=0.$ \end{proof}

 To proceed, let us use the following notation.

 \begin{notation}\label{n.fix}
In the setting of Proposition \ref{p.TZ},
 let $\Phi$ be the root system of $\fg$ with respect to $\ft$ and let $\fg_{\alpha} \subset \fg$  be the root space corresponding to a root $ \alpha \in \Phi.$ For a $\ft$-stable subspace $W \subset \fg$, we have the decomposition analogous to (\ref{e.decompo}) $$W = (W \cap \ft) \oplus (\oplus_{\alpha \in \Phi_W} \fg_{\alpha})$$ with $ \Phi_W := \{ \alpha \in \Phi \mid \fg_{\alpha} \subset W \}. $ From above, we have $$\Phi_{\fh} \subset \Phi_{\fn}, \ \Phi_{\fc_{\fg}(s)} \subset \Phi_{\fn}, \ \Phi = \Phi_{\fn} \sqcup \Phi_{\fm} \mbox{ and } \Phi_{T_s \widehat{Z}} = \Phi_{\fh} \sqcup \Phi_{\fm}.$$ \end{notation}

\begin{proposition}\label{p.nhc}
In the setting of Proposition \ref{p.TZ}, we have the following.
\begin{itemize}
\item[(i)] $\fn = \fh + \fc_{\fg}(s).$
\item[(ii)] $\Phi = \Phi_{\fm} \sqcup (\Phi_{\fh} \cup \Phi_{\fc_{\fg}(s)}).$
\item[(iii)] $\fc_{\fg}(\fh_0) \subset \fc_{\fg}(s) \subset \fh + \fc_{\fg}(\fh_0).$
\item[(iv)] $\fn = \fh + \fc_{\fg}(\fh_0).$
\item[(v)] If $\alpha \in \Phi_{\fm}$, then $[\fg_{\alpha},s] = \fg_{\alpha} \neq 0$ and consequently, the subspace $$\fa_{\alpha} := \{ a \in \fh_0 \mid [\fg_{\alpha}, a] =0\}$$ has codimension 1 in $\fh_0.$
\item[(vi)] If $\alpha \in \Phi$ satisfies $[\fg_{\alpha}, s] \neq 0$, then $\alpha \in \Phi_{\fm} \sqcup \Phi_{\fh}.$\end{itemize} \end{proposition}

\begin{proof}
In (i), the inclusion $\fh + \fc_{\fg}(s) \subset \fn$ is immediate. For any  $\alpha \in \Phi_{\fn} \setminus \Phi_{\fc_{\fg}(s)}$, we have $[\fg_{\alpha}, s] \neq 0$, which  implies that $$\fg_{\alpha} = [\fg_{\alpha}, s] \subset [\fn, \fh] \subset \fh,$$  hence $\alpha \in \Phi_{\fh}$. It follows that $\Phi_{\fn} \subset \Phi_{\fh}\cup \Phi_{\fc_{\fg}(s)}$ and $$ \fn \subset \ft + \fh + \fc_{\fg}(s) = \fh + \fc_{\fg}(s),$$ proving the reverse inclusion.

By (i), we have $\fg = \fn \oplus \fm = (\fh + \fc_{\fg}(s)) \oplus \fm$, which implies (ii).

In (iii), the first inclusion is immediate from $s \in \fh_0 \subset \ft$. It remains to show $\Phi_{\fc_{\fg}(s)} \subset \Phi_{\fc_{\fg}(\fh_0)} \cup \Phi_{\fh}$, or equivalently, any $\alpha \in \Phi_{\fc_{\fg}(s)} \setminus \Phi_{\fc_{\fg}(\fh_0)}$ belongs to $\Phi_{\fh}.$
From  $$0 = [\fg_{\alpha}, s] \mbox{ and }0 \neq [\fg_{\alpha}, \fh_0] \subset [\fg_{\alpha}, \ft] \subset \fg_{\alpha} ,$$ we can choose $w \in \fh_0 \setminus \C s$ and $0 \neq x_{\alpha} \in \fg_{\alpha}$ such that $[x_{\alpha}, w] = x_{\alpha}$. Define $s_{\epsilon} := s + \epsilon w \in \fh_0$ for $\epsilon \in \C$. By Corollary \ref{c.TZ}, we have for a general $\epsilon \in \C$, \begin{eqnarray*} T_{s_{\epsilon}} \widehat{Z} &=& \fh + [\fg, s_{\epsilon}] \\ &=& \fh + [\fn, s_{\epsilon}] + [\fm, s_{\epsilon}] \\ & \subset & \fh + \fm = T_s \widehat{Z}, \end{eqnarray*} using $[\fn, s_{\epsilon}] \subset [\fn, \fh] \subset \fh$ and $[\fm, s_{\epsilon}] \subset [\fm, \ft] \subset \fm$. But $\dim T_{s_{\epsilon}} \widehat{Z} = \dim T_s \widehat{Z}$ for a general $\epsilon$. Thus \begin{equation}\label{e.T=T} \fh + [\fg, s_{\epsilon}] = T_{s_{\epsilon}} \widehat{Z} = T_s \widehat{Z} = \fh + \fm.\end{equation} From $\alpha \in \Phi_{\fc_{\fg}(s)}$ and $x_{\alpha} \in \fc_{\fg}(s)$, we have $$ x_{\alpha} = [x_{\alpha}, w] = [x_{\alpha}, \frac{1}{\epsilon}(s_{\epsilon} - s)]  = \frac{1}{\epsilon} [ x_{\alpha}, s_{\epsilon}] \in [\fg, s_{\epsilon}].$$ By (\ref{e.T=T}), the last term lies in $\fh \oplus \fm$. Thus $\alpha \in (\Phi_{\fh} \sqcup \Phi_{\fm}) \cap \Phi_{\fc_{\fg}(s)} \subset \Phi_{\fh}$, completing the proof of (iii).

 (iv) follows from (i) and (iii).

Since $\fc_{\fg}(s) \subset \fn$ from Proposition \ref{p.tstable} (i), we have $\Phi_{\fm} \cap \Phi_{\fc_{\fg}(s)} = \emptyset.$ This implies (v).

(vi) is immediate from (ii).  \end{proof}

\begin{proposition}\label{p.nGauss}
In Set-up \ref{setup},  we have the Gauss map $Z \dasharrow {\rm Gr}(p+1, \fg)$ to the Grassmannian  of $(p+1)$-dimensional subspaces in $\fg$ for $p= \dim Z.$ Let $R \subset {\rm Gr}(p+1, \fg)$ be the proper image  of the Gauss map.  Assume that    the morphism $\gamma: Z^o \to R^o$ is the restriction of the Gauss map to suitable Zariski-open subsets $Z^o \subset Z$ and $R^o \subset R$.    Then in the notation of Proposition \ref{p.TZ},
\begin{itemize} \item[(i)]  $[\fg, \fh] \subset \fh \oplus \fm$; and \item[(ii)]  $\fn = \fh + \fc_{\fg}(\fh) + \ft. $ \end{itemize} \end{proposition}

\begin{proof}
For any general $v \in \fh$, we have $$\fh + [ \fg, v] = T_v \widehat{Z} = T_s \widehat{Z} = \fh \oplus \fm,$$ where the first equality is from Corollary \ref{c.TZ}, the second equality is from the assumption that $\gamma$ is the restriction of the Gauss map and the third equality is from Proposition \ref{p.TZ}. It follows that $[\fg, v] \subset \fh \oplus \fm$ for a general $v \in \fh$, proving (i).

In (ii), the inclusion $\fh + \fc_{\fg}(\fh) + \ft \subset \fh + \fc_{\fg}(\fh_0) = \fn$ is clear from $\ft + \fc_{\fg}(\fh) \subset \fc_{\fg}(\fh_0)$ and Proposition \ref{p.nhc} (iv). It remains to show $\Phi_{\fc_{\fg}(\fh_0)} \subset \Phi_{\fh} \cup \Phi_{\fc_{\fg}(\fh)}.$

Assume the contrary, that there exists $\alpha \in \Phi$ such that $$[\fg_{\alpha}, \fh_0] =0, \ [\fg_{\alpha}, \fh] \neq 0 \mbox{ and } \alpha \not\in \Phi_{\fh}.$$ From $\fh = \fh_0 \oplus (\oplus_{\delta \in \Phi_{\fh}} \fg_{\delta})$, we can choose $0 \neq x_\alpha \in \fg_{\alpha}$ and $0 \neq x_{\beta} \in \fg_{\beta}, \beta \in \Phi_{\fh}$ such that $ [x_{\alpha}, x_{\beta}] \neq 0. $ Note that \begin{equation}\label{e.23} [x_{\alpha}, x_{\beta}] \in [\fc_{\fg}(\fh_0), \fh] \subset \fh \end{equation} from Proposition \ref{p.nhc} (iv).

We claim that $\alpha \neq - \beta.$ Otherwise, setting $$t_{\alpha} := [x_{\alpha}, x_{\beta}] = [x_{\alpha}, x_{-\alpha}]  \in \ft,$$ we have the subalgebra in $\fg$  isomorphic to $\fsl_2(\C)$, spanned by $x_{\alpha}, x_{- \alpha}$ and $ t_{\alpha}$. In particular, we have $[t_{\alpha}, x_{\alpha}] \neq 0.$ But $t_{\alpha} \in \ft \cap \fh = \fh_0$ by (\ref{e.23}),  which yields $[t_{\alpha}, x_{\alpha}] \in [\fh_0, \fc_{\fg}(\fh_0)] =0$, a contradiction. This verifies the claim.

By the above claim and $[x_{\alpha}, x_{\beta}] \neq 0,$ we have $\alpha + \beta \in \Phi$. Then (\ref{e.23}) implies $\alpha + \beta \in \Phi_{\fh}$. Thus $$\fg_{\alpha} = [ \fg_{-\beta}, \fg_{\alpha + \beta}] \subset [\fg, \fh] \subset \fh \oplus \fm$$ by (i). Combining it with $\fg_{\alpha} \subset \fc_{\fg}(\fh_0) \subset \fn$ from Proposition \ref{p.nhc} (iv), we obtain $$\fg_{\alpha} \subset \fn \cap (\fh \oplus \fm) \subset \fn \oplus \fm = \fg.$$  This implies $\fg_{\alpha} \subset \fh$, a contradiction to our choice $\alpha \not\in \Phi_{\fh}.$
 \end{proof}

\begin{proposition}\label{p.dimh01}
In the setting of Proposition \ref{p.nGauss}, assume that $\fh_0 =\C s$. Then there exists a root $\alpha \in \Phi$ such that $$\fh = \fsl_2(\alpha) = \fg_{\alpha} + \ft_{\alpha} + \fg_{-\alpha} \mbox{ and } \fh_0 = \ft_{\alpha}.$$ \end{proposition}

\begin{proof}
As $\fh$ is $\ft$-stable from Proposition \ref{p.tstable} (ii), we can write $$\fh = \fh_0 \oplus \bigoplus_{\alpha \in \Phi_{\fh}} \fg_{\alpha}.$$ From the assumption $\dim \fh_0 =1$ and $\dim \fh \geq 2$ in Set-up \ref{setup} (3), we see $\Phi_{\fh} \neq \emptyset.$
We claim  \begin{equation}\label{e.beta} \fh_0 = \C s = \ft_{\beta} \mbox{ for any } \beta \in \Phi_{\fh}. \end{equation}
From Proposition \ref{p.TZ} and  Proposition \ref{p.nGauss} (i), we have $$0 \neq \ft_{\beta} = [\fg_{-\beta}, \fg_{\beta}] \ \subset  \ [\fg, \fh] \cap \ft \ \subset \ (\fh + \fm ) \cap \ft \  = \ \fh_0 = \C s.$$ This proves (\ref{e.beta}).

For any $\alpha, \alpha' \in \Phi_{\fh},$ (\ref{e.beta}) says $\ft_{\alpha} = \ft_{\alpha'}$, which implies $\alpha' = \alpha$ or $\alpha' = - \alpha$. Thus to prove Proposition \ref{p.dimh01}, we may assume that $\fh = \fg_{\alpha} + \ft_{\alpha}$ for some $\alpha \in \Phi$ and derive a contradiction.

The two subalgebras $\fh = \fg_{\alpha} + \ft_{\alpha}$ and $\fh':= \fg_{-\alpha} + \ft_{\alpha}$ are opposite Borel subalgebras of the simple Lie algebra $\fsl_2(\alpha) = \fg_{\alpha} + \ft_{\alpha} + \fg_{-\alpha}$. Thus there exists $g \in \exp (\fsl_2(\alpha)) \subset G$ such that $${\rm Ad}_g (\ft_{\alpha}) = \ft_{\alpha}, \ {\rm Ad}_g(\fg_{\alpha}) = \fg_{-\alpha}, \mbox{ and } {\rm Ad}_g(\fg_{-\alpha}) = \fg_{\alpha}.$$ This means that $g \cdot [s] = [s]$ and $g \cdot \BP \fh = \BP \fh'$. But $\gamma$ is $G$-equivariant, which implies $$g \cdot r = g \cdot \gamma ([s]) = \gamma( g \cdot [s]) = \gamma ([s]) =r.$$ Consequently,  the closure $\BP \fh'$ of $\gamma^{-1}(g \cdot r)$ agrees with the closure $\BP \fh$ of $\gamma^{-1}(r)$  for $r := \gamma ([s]) \in R^o$. This leads to the contradiction $\fh = \fh'$.
\end{proof}

For the next proposition, we need the following general fact.

\begin{lemma}\label{l.SL2}
 Let $G$ be a simple algebraic group. Fix a choice of a system of simple roots  of the Lie algebra $\fg$ of $G$.  For a simple root $\alpha$,  let $G(\alpha) \subset G$ be the connected subgroup whose Lie algebra is the subalgebra  $\fg_{\alpha} + \ft_{\alpha} + \fg_{-\alpha}$ of $\fg$ spanned by the root spaces of $\pm \alpha$ and $\ft_{\alpha} := [\fg_{\alpha}, \fg_{-\alpha}]$. Assume that $(\fg, \alpha)$ is different from the following two cases: \begin{itemize} \item   $\fg \cong \fsl_2(\C)$, or \item   $\fg$ is of type $B_r$ and $\alpha= \alpha_r$ is the unique short simple root. \end{itemize} Then  $G(\alpha)$ is isomorphic to ${\rm SL}_2(\C)$.   \end{lemma}

 \begin{proof}
 The assumption implies that there exists another simple root $\beta $ such that $\langle \beta, \alpha\rangle = -1$ or $-3$. If $\langle \beta, \alpha \rangle =-1$ (resp. $-3$), then there are $2$ roots (resp. $4$ roots) in  the $\alpha$-string through $\beta$ and their root subspaces span  an irreducible $G(\alpha)$-submodule in $\fg$ of   dimension $2$  (resp. $4$).   The three-dimensional group $G(\alpha)$ is isomorphic to either ${\rm SL}_2(\C)$ or ${\rm PSL}_2(\C)$.  Since ${\rm PSL}_2(\C)$ cannot have an even-dimensional irreducible representation, the group $G(\alpha)$ must be isomorphic to ${\rm SL}_2(\C)$. \end{proof}

\begin{proposition}\label{p.SL2}
In the setting of Proposition \ref{p.dimh01}, the connected subgroup $G(\alpha) \subset G$ corresponding to the Lie subalgebra $\fsl_2(\alpha) \subset \fg$ is isomorphic to ${\rm SL}_2(\C).$ \end{proposition}

\begin{proof} We fix a system of simple roots including $\alpha.$
By Lemma \ref{l.SL2}, we need to exclude the following two possibilities.
\begin{itemize}  \item   $\fg \cong \fsl_2(\C)$, or \item   $\fg$ is of type $B_r$ and $\alpha= \alpha_r$ is the unique short simple root. \end{itemize}
If $\fg \cong \fsl_2(\C)$, Proposition  \ref{p.dimh01} implies that $\fh \cong \fg$. Then $Z = \BP \fg$ and the Gauss map $\gamma$ must be a constant, a contradiction.

If $\fg$ is of type $B_r$ and $\alpha = \alpha_r$ is the unique short simple root, set $\beta:= \alpha_{r-1} + \alpha_r.$ Then $\langle \beta, \alpha \rangle =0$ implies $\fg_{\beta} \subset \fc_{\fg}(\ft_{\alpha}) = \fc_{\fg}(\fh_0).$ By Proposition \ref{p.nGauss},
$$\fg_{\beta} \ \subset \ \fc_{\fg}(\fh_0) \ \subset \ \fn \ = \ \fh + \fc_{\fg}(\fh) + \ft.$$
It follows that $$ \alpha_{r-1} + \alpha_r = \beta \ \in \  \Phi_{\fh} \sqcup \Phi_{\fc_{\fg}(\fh)} = \{ \pm \alpha\} \sqcup \Phi_{\fc_{\fg}(\fh)},$$ which implies $\beta \in \Phi_{\fc_{\fg}(\fh)}$ and $\fg_{\beta} \subset \fc_{\fg}(\fh)$. Then $$0 \neq \fg_{\alpha_{r-1}} = [\fg_{\alpha_{r-1} + \alpha_r}, \fg_{-\alpha_r}] \subset [\fg_{\beta}, \fh] =0,$$ a contradiction. \end{proof}

The assumption $\dim \fh_0 =1$ in Proposition \ref{p.dimh01} arises from the following geometric condition.

\begin{proposition}\label{p.sn}
In the setting of Proposition \ref{p.nGauss}, assume that the normalization of the variety $Z$ is nonsingular. Then $\dim \fh_0 =1.$ \end{proposition}

\begin{proof}
Assuming $\dim \fh_0 \geq 2$, we derive a contradiction as follows.

Let $\nu: Y \to Z$ be the normalization of $Z$ with $Y$ nonsingular by our assumption.  Let $F \subset Y$ be the proper transform of $\BP \fh \subset Z$.  The restriction $\nu|_F : F \to \BP \fh$ is a finite birational morphism,  because we are assuming that $\BP \fh$ is the closure of a general fiber of  $\gamma$. It follows that $\nu|_F$ must be biregular and $F$ is isomorphic to the projective space $\BP \fh$.

From the hypothesis $\dim \fh_0 \geq 2$,  we can choose a general $2$-dimensional subspace $\fk \subset \fh_0$ containing $s$.
Let $\ell\subset F$ be the inverse image $(\nu|_F)^{-1}(\BP \fk)$, which is a line on $F$.
We have a natural $G$-action on $Y$ lifting the $G$-action on $Z$ given by the adjoint action on $\fg$.
For each $\beta \in \Phi_{\fm}$, the intersection $\fk \cap \fa_{\beta}$ with the subspace $\fa_{\beta} \subset \fh_0$ of codimension 1 in Proposition \ref{p.nhc} (v) has dimension $1$. Hence    the subgroup $\exp(\fg_{\beta}) \subset G$ transforms $\ell = \ell_0$ in a one-dimensional family of rational curves $\{\ell_t \subset Y \mid t \in \C\}$ such that \begin{itemize} \item[(1)]   the point $P_{\beta} \in \ell$ corresponding to $\BP (\fa_{\beta} \cap \fk) \in \BP \fk$ is fixed by the action of $\exp(\fg_{\beta})$, namely, it is contained in $\ell_t$ for all $t \in \C$; and \item[(2)] the tangent space of the $\exp(\fg_{\beta})$-orbit through the point  $y \in \ell$ corresponding to $s$ spans the 1-dimensional subspace $[\fg_{\beta}, s] =\fg_{\beta}$ from Proposition \ref{p.nhc} (v). \end{itemize}
The infinitesimal deformation of the family $\{\ell_t \subset Y \mid t \in \C\}$ gives a holomorphic section $\sigma_{\beta}$ of the normal bundle $N_{\ell/Y}$ such that \begin{itemize} \item[(a)]
$\sigma_{\beta}(P_{\beta}) =0$ and \item[(b)]  the differential ${\rm d}_y \nu: T_y Y \to T_z Z, z = \nu(y), s \in \widehat{z},$ induces an isomorphism $$N_{\ell/Y, y} \ \stackrel{\approx}{\longrightarrow} \  N_{\BP \fk/Z, z} = \Hom(\C s,  T_s \widehat {Z}/\fk)$$ that sends $ \C \sigma_{\beta}(y)$ to the one-dimensional subspace  $$\Hom(\C s, \fg_{\beta}) \ \subset \ \Hom(\C s, \fm) \ \subset \ \Hom(\C s, T_s \widehat{Z}/\fk),$$ where  we regard $\fm$ as a subspace of $T_s \widehat{Z}/\fk$ using $$(T_s \widehat{Z} \cap \fk)\cap (T_s \widehat{Z} \cap \fm) = 0 \mbox{  from }  \fk \cap \fm \subset \fh_0 \cap \fm =0.$$ \end{itemize}   Since the normal bundle $N_{\ell/F}$ is isomorphic to $\sO(1)^{\oplus (\dim F -1)}$ and the sections  $\sigma_{\beta}, \beta \in \Phi_{\fm}$ of the normal bundle $N_{\ell/Y}$ vanishing at $P_{\beta} \in \ell$ are transversal to $F$ at the point $y$, the existence of sections $\sigma_{\beta}$ satisfying (a) and (b) for all $\beta \in \Phi_{\fm}$ implies that the normal bundle $N_{\ell/Y}$ is isomorphic to $\sO(1)^{\oplus (\dim Y -1)}$.

Since the line bundle $\nu^* \sO_{\BP \fg}(1)$ is ample on $Y$ and has intersection number 1 with $\ell \subset Y$, we see that $\ell$ is a minimal rational curve on the nonsingular projective variety $Y$ with the normal bundle isomorphic to  $\sO(1)^{\oplus (\dim Y -1)}$. It follows that $Y$ is isomorphic to $\BP^{\dim Z}$ and the line bundle $\nu^* \sO_{\BP \fg}(1)$  is the hyperplane line bundle on $Y$ (for example, by \cite[Corollary 0.3]{CMS}). It follows that $Z = \nu(Y) \subset \BP \fg$ must be a linear subspace and the Gauss map $\gamma$ is constant, a contradiction to the assumption that $R$ is positive-dimensional.  \end{proof}

Now we are ready to prove Lemma \ref{l.technical}.

\begin{proof}[Proof of Lemma \ref{l.technical}]
Under the assumptions of Lemma \ref{l.technical}, we are in the setting of Proposition \ref{p.sn}.  Thus we have $\dim \fh_0 =1$ and  can use Propositions \ref{p.dimh01} and \ref{p.SL2} to conclude that  $\exp(\fh) \subset G$ is isomorphic to ${\rm SL}_2(\C).$ \end{proof}

\medskip
{\bf Acknowledgment} We would like to thank Baohua Fu and Yoshinori Namikawa for valuable discussions on nilpotent orbits.

\medskip
Jun-Muk Hwang (jmhwang@ibs.re.kr)
Center for Complex Geometry,
Institute for Basic Science (IBS),
Daejeon 34126, Republic of Korea

\medskip
Qifeng Li (qifengli@sdu.edu.cn)
School of Mathematics,
Shandong University,
Jinan 250100, China
\end{document}